\documentclass[a4paper,10pt]{article}
\usepackage[dvips]{graphicx}
\usepackage{float}
\usepackage{amsmath} 
\usepackage{amssymb}
\usepackage{amsfonts}
\usepackage{theorem}

\numberwithin{equation}{section}

\newtheorem{theo}{Theorem}[section]
\newtheorem{lem}{Lemma}[section]
\newtheorem{cor}{Corollary}[section]
\newtheorem{pro}{Proposition}[section]

{
\theoremstyle{break}
\theorembodyfont{\rmfamily}
}

\renewcommand{\thefootnote}{\fnsymbol{footnote}}

\newcommand{\bX}{\mathbf{X}}
\newcommand{\bx}{\mathbf{x}}

\newif\ifnotes
\notestrue
\newcounter{mnotecount}[section]

\begin{document}

\begin{center}

{\sc \Large PCA-Kernel Estimation\\
\vspace{0.5cm}}

G\'erard Biau
\vspace{0.2cm}\\
LSTA \& LPMA\\
Universit\'e Pierre et Marie Curie -- Paris VI\\
Bo\^{\i}te 158,   175 rue du Chevaleret\\
75013 Paris, France

\vspace{0.2cm}
DMA\\
Ecole Normale Sup\'erieure\\
45 rue d'Ulm\\
75230 Paris Cedex 05, France\\
\smallskip
\texttt{gerard.biau@upmc.fr}
\vspace{0.5cm}

Andr\'e Mas\footnote{Corresponding
author.} 
\vspace{0.2cm}\\
Institut de Math\'ematiques et de Mod\'elisation de
Montpellier\\ 
UMR CNRS 5149, Equipe de Probabilit\'es et Statistique\\
Universit\'e Montpellier II, CC 051\\ 
Place Eug\`ene Bataillon, 34095 Montpellier Cedex 5, France\\
\smallskip
\texttt{mas@math.univ-montp2.fr}

\end{center}

\begin{abstract}

\medskip 
\noindent {\rm 
Many statistical estimation techniques for high-dimensional or functional data are based on a preliminary dimension reduction step, which consists in projecting the sample $\bX_1, \hdots, \bX_n$ onto the first $D$ eigenvectors of the Principal Component Analysis (PCA) associated with the empirical projector $\hat \Pi_D$.  Classical nonparametric inference methods such as kernel density estimation or kernel regression analysis are then performed in the (usually small) $D$-di\-men\-sio\-nal space. However, the mathematical analysis of this data-driven dimension reduction scheme raises technical problems, due to the fact that the random variables of the projected sample $( \hat \Pi_D\bX_1,\hdots, \hat \Pi_D\bX_n )$ are no more independent. As a reference for further studies, we offer in this paper several results showing the asymptotic equivalencies between important kernel-related quantities based on the empirical projector and its theoretical counterpart. As an illustration, we provide an in-depth analysis of the nonparametric kernel regression case.\\

\medskip
 
}

\end{abstract}

\renewcommand{\thefootnote}{\arabic{footnote}}

\setcounter{footnote}{0}
\section{Introduction}
Nonparametric curve estimation provides a useful tool for exploring and understanding the structure of a data set, especially when parametric models are inappropriate. A large amount of progress has been made in the 90's in both the design and the study of inferential aspects of nonparametric estimates. There are too many references to be included here, but the monographs of Silverman \cite{Silverman}, Scott \cite{Scott}, Simonoff \cite{Simonoff} and Gy\"orfi et al. \cite{Laciregressionbook} will provide the reader with good introductions to the general subject area. \\

Among all the nonparametric methods which have been proposed so far, kernel estimation has gained favor from many data analysts, probably because of its simplicity to implement and good statistical properties---see for example Simonoff \cite{Simonoff} for a variety of real data examples which illustrate the power of the approach. Kernel estimates were originally studied in density estimation by Rosenblatt \cite{Rosenblatt} and Parzen \cite{Parzen}, and were latter introduced in regression estimation by Nadaraya \cite{Nadaraya1, Nadaraya2} and Watson \cite{Watson}. A compilation of the mathematical properties of kernel estimates can be found in  Prakasa Rao \cite{Rao} (for density estimation), Gy\"orfi et al. \cite{Laciregressionbook} (for regression) and Devroye et al. \cite{DGL} (for classification and pattern recognition). To date, most of the results pertaining to kernel estimation have
been reported in the finite-dimensional case, where it is assumed that the observation space is
the standard Euclidean space $\mathbb R^d$. However, in an increasing number of practical applications, input data items
are in the form of random functions (speech recordings, multiple time series,
images...) rather than standard vectors, and this casts the problem
into the general class of functional data analysis. Motivated by this broad range of potential applications, Ferraty and Vieu describe in \cite{FV} a possible route to extend kernel estimation to potentially infinite-dimensional spaces.\\

On the other hand, it has become increasingly clear over the years that the performances of kernel estimates deteriorate as the dimension of the problem increases. The reason for this is that, in high dimensions, local neighborhoods tend to be empty of sample observations unless the sample size is very large. Thus, in kernel estimation,  there will be no local averages to take unless the bandwidth is very large. This general problem was termed the curse of dimensionality (Bellman \cite{Bellman}) and, in fact, practical and theoretical arguments suggest that kernel estimation beyond 5 dimensions is fruitless. The paper by Scott and Wand \cite{SW} gives a good account on the feasibility and difficulties of high-dimensional estimation, with examples and computations.\\

In order to circumvent the high-dimension difficulty and make kernel estimation simpler, a wide range of techniques have been developed. One of the most common approaches is a two-stage strategy: first reduce the dimension of the data and then perform---density or regression---kernel estimation. With this respect, a natural way to reduce dimension is to extract the largest $D$ principal component axes (with $D$ chosen to account for most of the variation in the data), and then operate in this $D$-dimensional space, thereby improving the ability to discover interesting structures (Jee \cite{Jee}, Friedman \cite{Friedman} and Scott \cite{Scott}, Chapter 7). To illustrate more formally this mechanism, let $(\mathcal F, \langle.,.\rangle,\|.\|)$ be a (typically high or infinite-dimensional) separable Hilbert space, and consider for example the regression problem, where we observe a set $\mathcal D_n=\{(\bX_1,Y_1), \hdots, (\bX_n,Y_n)\}$ of independent $\mathcal F\times \mathbb R$-valued random variables with the same distribution as a generic pair $(\bX, Y)$ satisfying $\mathbb E |Y|<\infty$. The goal is to estimate the regression function $r(\bx)=\mathbb E[Y|\bX=\bx]$ using the data $\mathcal D_n$. The kernel estimate of the function $r$ takes the form
$$r_n(\bx)=\frac{\sum_{i=1}^n Y_iK\left(\frac{\|\bx-\bX_i\|}{h_n} \right)}{\sum_{i=1}^nK\left(\frac{\|\bx-\bX_i\|}{h_n} \right)}$$
if the denominator is nonzero, and $0$ otherwise. Here the bandwidth $h_n>0$ depends only on the sample size $n$, and the function $K:[0,\infty) \to [0,\infty)$ is called a kernel. Usually, $K(v)$ is ``large'' if $v$ is ``small'', and the kernel estimate is therefore a local averaging estimate. Typical choices for $K$ are the naive kernel $K(v)=\mathbf 1_{[0,1]}(v)$, the Epanechnikov kernel $K(v)=(1-v^2)_+$, and the Gaussian kernel $K(v)= \exp(-v^2/2)$.  \\

As explained earlier, the estimate $r_n$ is prone to the curse of dimensionality,  and the strategy advocated here is to first reduce the ambient dimension by the use of Principal Component Analysis (PCA, see for example Dauxois et al. \cite{Dauxois} and Jolliffe \cite{Jolliffe}). More precisely, assume without loss of generality that $\mathbb E \bX=0$, $\mathbb E\|\bX\|^2<\infty$, and let $\Gamma(.)  =\mathbb{E}[  \langle \bX,\cdot \rangle \bX] $ be the covariance operator of $\bX$ and  $\Pi_{D}$ be the orthogonal projector on the collection of the first $D$ eigenvectors $\{\mathbf e_1, \hdots, \mathbf e_D\}$ of $\Gamma$ associated with the first $D$ eigenvalues $\lambda_1 \geq \lambda_2Ê\geq \hdots\geq \lambda_D\geq 0$. In the sequel we will assume as well that the distribution of $\bX$ is nonatomic. \\

In this context, the PCA-kernel regression estimate reads
$$r_n^D(\bx)=\frac{\sum_{i=1}^n Y_iK\left(\frac{\|\Pi_D(\bx-\bX_i)\|}{h_n} \right)}{\sum_{i=1}^nK\left(\frac{\|\Pi_D(\bx-\bX_i)\|}{h_n} \right)}.$$
The hope here is that the most informative part of the distribution of $\bX$ should be preserved by projecting the observations on the first $D$ principal component axes, so that the estimate should still do a good job at estimating $r$ while performing in a  reduced-dimensional space. Alas, on the practical side, the smoother $r_n^D$ is useless since the distribution of $\bX$ (and thus, the projector $\Pi_D$) is usually unknown, making of $r_n^D$ what is called a ``pseudo-estimate''. However, the covariance operator $\Gamma$ can be approximated by its empirical version 
\begin{equation}
\label{empgamma}
\Gamma_n(.) =\frac{1}{n}\sum_{i=1}^n \langle \bX_i,\cdot \rangle \bX_i,
\end{equation}
and $\Pi_D$ is in turn approximated by the empirical orthogonal projector $\hat \Pi_D$ on the (empirical) eigenvalues $\{\hat {\mathbf e}_1, \hdots,\hat{\mathbf e}_D\}$ of $\Gamma_n$. Thus, the operational version $\hat r_n^D$ of the pseudo-estimate $r_n^D$ takes the form
$$\hat r_n^D(\bx)=\frac{\sum_{i=1}^n Y_iK\left(\frac{\|\hat \Pi_D(\bx-\bX_i)\|}{h_n} \right)}{\sum_{i=1}^nK\left(\frac{\|\hat \Pi_D(\bx-\bX_i)\|}{h_n} \right)}.$$
Unfortunately, from a mathematical point of view, computations involving the numerator
or the denominator of the estimate $\hat r_n^D$ are difficult, since the random variables $(K({\|\hat \Pi_D(\bx-\bX_i)\|}/{h_n})) _{1\leq i\leq n}$ are identically distributed but clearly {\it not} independent. Besides, due to nonlinearity, the distribution of
  $\|\hat \Pi_D(\bx-\bX_i)\|$ is usually inaccessible, even when the $\bX_{i}$'s have known and simple distributions. In short, this makes any theoretical calculation impossible, and it essentially explains why so few theoretical results have been reported so far on the statistical properties of the estimate $\hat r_n^D$, despite its wide use. On the other hand, we note that the random variables $(K({\|\Pi_D(\bx-\bX_i)\|}/{h_n})) _{1\leq i\leq n}$ are independent and identically distributed. Therefore, the pseudo-estimate $r_n^D$ is amenable to mathematical analysis, and fundamental asymptotic theorems such that the law of large numbers and the central limit theorem may be applied.\\
  
In the present contribution, we prove that $\hat r_n^D$ and $r_n^D$ have the same asymptotic behavior and show that nothing is lost in terms of rates of convergence when replacing $\hat r_n^D$ by $r_n^D$ (Section 4). In fact, taking a more general view, we offer in Section 3 a thorough asymptotic comparison of the partial sums
$$S_n(\bx)=\sum_{i=1}^n K\left(  \frac{\left\| \Pi_{D}\left(  \bx-\bX_{i}\right)  \right\|}{
h_n}\right) \quad \mbox{and}\quad \hat S_n(\bx)=\sum_{i=1}^n K\left(  \frac{\left\| \hat \Pi_{D}\left(  \bx-\bX_{i}\right)  \right\|}{
h_n}\right)$$
with important consequences in kernel density estimation. As an appetizer, we will first carry out in Section 2 a preliminary analysis of the asymptotic proximity between the projection operators $\hat \Pi_D$ and $\Pi_D$. Our approach will strongly rely on the representation of the operators by Cauchy integrals, through what is classically known in analysis as perturbation method. For the sake of clarity, proofs of the most technical results are postponed to Section 5.
\section{Asymptotics for PCA projectors}
Here and in the sequel, we let $(\mathcal F, \langle.,.\rangle,\|.\|)$ be a separable Hilbert space and $\bX_1, \hdots, \bX_n$ be independent random variables, distributed as a generic non\-ato\-mic and centered random $\bX$ satisfying $\mathbb E\|\bX\|^2<\infty$. Denoting by $\Gamma$ the covariance operator of $\bX$, we let $\Pi_D$ be the orthogonal projection operator on $\{\mathbf e_1, \hdots, \mathbf{e}_D\}$, the set of first $D$ eigenvectors of $\Gamma$ associated with the (nonnegative) eigenvalues $\{\lambda_1,\hdots,\lambda_D\}$ sorted by decreasing order. The empirical version $\Gamma_n$ of $\Gamma$ is defined in (\ref{empgamma}), and we denote by $\{\hat {\mathbf e}_1, \hdots, \hat {\mathbf e}_D\}$ and $\{\hat \lambda_1, \hdots, \hat \lambda_D\}$ the associated empirical eigenvector and (nonnegative) eigenvalue sets, respectively, based on the sample $\bX_1, \hdots, \bX_n$. To keep things simple, we will assume throughout that the projection dimension $D$ is fixed and independent of the observations (for data-dependent methods regarding the choice of $D$, see for example Jolliffe \cite{Jolliffe}). Besides, and without loss of generality, it will also be assumed that $\lambda_1>\hdots > \lambda_{D+1}$. This assumption may be removed at the expense of more tedious calculations taking into account the dimension of the eigenspaces (see for instance \cite{Mas} for a generic method). \\
 
\begin{center}
\begin{figure}[!h]
\includegraphics*[width=10cm,height=6cm]{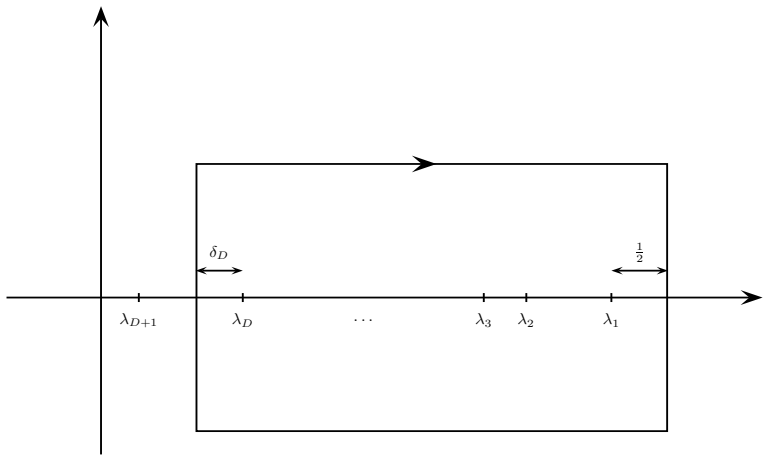}
\caption{Oriented rectangle-contour $\mathcal{C}_{D}$, with a right vertex intercepting the real line at $x=\lambda_1+1/2$ and a left vertex passing through $x=\lambda_{D}-\delta_D$, $\delta_{D}=(\lambda_{D}-\lambda_{D+1})/2.$\label{figure1}}
\end{figure}
\end{center}

The aim of this section is to derive new asymptotic results regarding the empirical projector $\hat \Pi_D$ on $\{\hat {\mathbf e}_1, \hdots, \hat {\mathbf e}_D\}$ as the sample size $n$ grows to infinity. Let us first recall some elementary facts from complex analysis. The eigenvalues $\{\lambda_1, \hdots, \lambda_D\}$ are nonnegative real numbers, but we may view them as points in the complex plane $\mathbb C$. Denote by $\mathcal{C}$ a closed oriented contour in $\mathbb{C}$, that is a closed curve (for instance, the boundary of a rectangle) endowed with a circulation. Suppose first that $\mathcal{C=}{\mathcal{C}}_{1}$ contains ${\lambda}_{1}$ only. Then, the so-called formula of residues (Rudin \cite{Rudin}) asserts that
\[
\int_{{\mathcal{C}}_{1}}\frac{\mbox{d}z}{z-{\lambda}_{1}}=1\quad \mbox{and} \quad \int_{{\mathcal{C}}_{1}}\frac{\mbox{d}z}{z-{\lambda}_{i}}=0 \quad \mbox{for } i\neq1.
\]
In fact, this formula may be generalized to functional calculus for
operators. We refer for instance to Dunford and Schwartz \cite{DS} or Gohberg et al. \cite{GGK} for exhaustive information about this theory, which allows to derive integration formulae for functions with operator values, such as 
\[
{\Pi}_{1}=\int_{{\mathcal{C}}_{1}}\left(  zI-\Gamma\right)  ^{-1}\mbox{d}z.
\]
Thus, in this formalism, the projector on ${\bf e}_1$ is explicitly written as a function of the covariance operator. Clearly, the same arguments allow to express the empirical projector $\hat \Pi_1$ as
$$\hat \Pi_{1}=\int_{\hat{\mathcal{C}}_{1}}\left(  zI-\Gamma_n\right)  ^{-1}\mbox{d}z,$$
where $\hat{\mathcal{C}}_{1}$ is a (random) contour which contains $\hat{\lambda}_{1}$ and no other
eigenvalue of $\Gamma_n$. These formulae generalize and, letting ${\mathcal C}_D$ (respectively $\hat{\mathcal C}_D$) be contours containing $\{\lambda_1, \hdots, \lambda_D\}$ (respectively $\{\hat{\lambda}_1, \hdots, \hat{\lambda}_D\}$) only, we may write 
\[
{\Pi}_{D}=\int_{{\mathcal{C}}_{D}}\left(  zI-\Gamma\right)  ^{-1}\mbox{d}z\quad \mbox{and}\quad \hat \Pi_{D}=\int_{\hat{\mathcal{C}}_{D}}\left(  zI-\Gamma_n
\right)  ^{-1}\mbox{d}z.
\]
The contours $\mathcal{C}_{D}$ may take different forms. However, to keep things simple, we let in the sequel $\mathcal{C}_{D}$ be the boundary of a rectangle as in Figure \ref{figure1}, with a right vertex intercepting the real line at $x=\lambda_1+1/2$ and a left vertex passing through $x=\lambda_{D}-\delta_D$, with 
$$\delta_{D}=\frac{\lambda_{D}-\lambda_{D+1}}{2}.$$
With a slight abuse of notation, we will also denote by $\mathcal C_D$ the corresponding rectangle.
Thus, with this choice, $\mathcal C_D$ contains $\{\lambda_1, \hdots, \lambda_D\}$ an no other eigenvalue. Lemma \ref{L0} below, which is proved in Section 5, shows that, asymptotically, this assertion is also true with $\{\hat{\lambda}_1, \hdots, \hat{\lambda}_D\}$ in place of $\{\lambda_1, \hdots, \lambda_D\}$. In the sequel, the letter $C$ will denote a positive constant, the value of
which may vary from line to line. Moreover, the notation $\left\| \cdot\right\| _{\infty}$ and $\left\| \cdot\right\| _{2}$ will stand for the classical operator and Hilbert-Schmidt norms, which are respectively defined by
$$\left\| T\right\| _{\infty}   =\sup_{\bx\in\mathcal{B}_{1}}\left\|
T\bx\right\| \quad \mbox{and} \quad  \left\| T\right\| _{2}^{2}   =\sum_{p=1}^{\infty}\left\|
T{\bf u}_{p}\right\| ^{2},$$
where $\mathcal{B}_{1}$ denotes the closed unit ball of $\mathcal F$ and $\left(
{\bf u}_{p}\right)  _{p \geq 1}$ a Hilbertian basis of $\mathcal F$. It is known (Dunford and Schwartz \cite{DS}) that the value of $\left\| T\right\| _{2}$ does not depend on the actual basis and that $\left\| \cdot\right\| _{\infty}\leq\left\| \cdot\right\| _{2}$. The Hilbert-Schmidt norm is of more generalized use, essentially because it yields simpler calculations than the sup-norm.\\

As promised, the next lemma ensures that the empirical eigenvalues are located in the rectangle $\mathcal{C}_{D}$ through an exponential concentration inequality.
\begin{lem}
\label{L0}For all $n \geq 1$, let the event
$$\mathcal{A}_{n}=\left\{  \hat{\lambda}_{i}\in
\mathcal{C}_{D},\,i=1, \hdots, D,\mbox{and }\hat{\lambda}_{D+1}\notin\mathcal{C}%
_{D}\right\}.  $$
There exists a positive constant $C$ such that
$$ \mathbb P(\mathcal A_n^c) =\mathcal O\left(\exp(-Cn)\right).$$
\end{lem}
Remark that the constants involved in the document depend on the actual dimension $D$ and their values increase as $D$ becomes large. To circumvent this difficulty, a possible approach is to let $D$ depend on $n$. This is beyond the scope of the present paper, and we refer to Cardot et al. \cite{CMS} for some perspectives in this direction.\\

We are now in a position to state the main result of the section. Theorem \ref{T1} below states the asymptotic proximity of the operators $\hat \Pi_D$ and $\Pi_D$, as $n$ becomes large, with respect to different proximity criteria. We will make repeated use of this result throughout the document. We believe however that it is interesting by itself. For a sequence of random variables $(Z_n)_{n\geq 1}$ and a positive sequence $(v_n)_{n \geq 1}$, notation $Z_n=\mathcal O(v_n)$ a.s. means that each random draw of $Z_n$ is $\mathcal O(v_n)$.
\begin{theo}
\label{T1}  The following three assertions are true for all $n \geq 1$:
\begin{itemize}
\item[$(i)$] There exists a positive constant $C$ such that, for all $\varepsilon >0$,
\begin{equation*}
\mathbb{P}\left(  \left\| \hat{\Pi}_{D}-\Pi_{D}\right\|_{\infty
}\geq \varepsilon\right)  = \mathcal O\left(\exp\left(  -Cn\varepsilon^{2}\right)\right).
\label{inegexp-proj}
\end{equation*}
\item[$(ii)$] One has
\begin{equation*}
\left\| \hat{\Pi}_{D}-\Pi_{D}\right\| _{\infty}  =\mathcal {O}\left(
\sqrt{\frac{\log n}{n}}\right)  \quad \mbox{a.s.}
\label{conv-proj-ps}
\end{equation*}
\item[$(iii)$] One has
\begin{equation*}
\mathbb{E}\left\| \hat{\Pi}_{D}-\Pi_{D}\right\| _{\infty}^{2}=\mathcal {O}\left(  \frac{1}{n}\right) .
\label{conv-proj-prob}%
\end{equation*}
\end{itemize}
\end{theo}
\noindent \textbf{Proof of Theorem \ref{T1}}\quad The proof will be based on arguments presented in Mas and Menneteau \cite{Mas}. Using the notation of Lemma \ref{L0}, we start from the decomposition%
\begin{equation}
\label{noel}
\hat{\Pi }_{D}-\Pi _{D}=\left( \hat{\Pi }_{D}-\Pi _{D}\right)\mathbf 1_{\mathcal{A}_{n}^{c}}+\left( \hat{\Pi }_{D}-\Pi _{D}\right)\mathbf 1_{\mathcal{A}_{n}}.
\end{equation}%
Consequently,
\begin{align*}
\mathbb{P} & \left(  \left\| \hat{\Pi}_{D}-\Pi_{D}\right\|_{\infty
}\geq \varepsilon\right)\\
&    \leq \mathbb{P}\left(  \left\| \hat{\Pi}_{D}-\Pi_{D}\right\|_{\infty
} \mathbf 1_{\mathcal{A}%
_{n}^{c}}\geq \varepsilon/2\right) + \mathbb{P}\left(  \left\| \hat{\Pi}_{D}-\Pi_{D}\right\|_{\infty
} \mathbf 1_{\mathcal{A}%
_{n}}\geq \varepsilon/2\right).
\end{align*}
Observing that
$$\left\|\hat{\Pi}_{D}-\Pi_{D} \right \| _{\infty } \mathbf 1_{\mathcal{A}%
_{n}^{c}}\leq 2\mathbf 1_{\mathcal{A}_{n}^{c}},$$
we conclude by Lemma \ref{L0} that
\begin{align}
\mathbb{P}\left( \left\|\hat{\Pi}_{D}-\Pi_{D} \right \| _{\infty } \mathbf 1_{\mathcal{A}_{n}^{c}} \geq \varepsilon/2\right) & \leq \mathbb{P}\left( \mathcal{A}_{n}^{c}\right) \notag\\
& =\mathcal O\left(\exp(-nC\varepsilon^2)\right). \label{ppl}
\end{align}
With respect to the second term in (\ref{noel}), write
\begin{align*}
\left(\hat{\Pi}_{D}-\Pi_{D}\right)\mathbf 1_{\mathcal A_n}  & =\mathbf 1_{\mathcal A_n}\int_{\mathcal{C}_{D}}\left[  \left(
zI-\Gamma_{n}\right)  ^{-1}-\left(  zI-\Gamma\right)  ^{-1}\right]  \mbox{d}z\\
& =\mathbf 1_{\mathcal A_n}\int_{\mathcal{C}_{D}}\left[  \left(  zI-\Gamma_{n}\right)  ^{-1}\left(
\Gamma_n-\Gamma\right)  \left(  zI-\Gamma\right)  ^{-1}\right]  \mbox{d}z.
\end{align*}
Let $\ell_{D}$ be the length of the contour $\mathcal{C}_{D}$. Using elementary properties of Riesz integrals (Gohberg et al. \cite{GGK}), we obtain
\begin{align*}
&\left\| \hat{\Pi}_{D}-\Pi_{D}\right\| _{\infty}\mathbf 1_{\mathcal A_n}\\
& \quad \leq\ell%
_{D}\left\| \Gamma_n-\Gamma\right\| _{\infty}\sup_{z\in\mathcal{C}%
_{D}}\left[  \left\| \left(  zI-\Gamma_{n}\right)  ^{-1}\right\|
_{\infty}\left\| \left(  zI-\Gamma\right)  ^{-1}\right\| _{\infty
}\right]\mathbf 1_{\mathcal A_n}.
\end{align*}
Observing that the eigenvalues of the symmetric operator $\left(  zI-\Gamma\right)
^{-1}$ are the $\left\{  \left(  z-\lambda_{i}\right)  ^{-1},i\in
\mathbb{N^{\star}}\right\}$, we see that $\left\| \left(  zI-\Gamma\right)  ^{-1}\right\|
_{\infty}=\mathcal{O}(\delta_{D})$. The same bound is valid taking $\Gamma_{n}$
instead of $\Gamma$, when $\mathcal{A}_{n}$ holds. In consequence, 
\begin{equation}
\left\| \hat{\Pi}_{D}-\Pi_{D}\right\| _{\infty}\mathbf 1_{\mathcal A_n}=\mathcal {O}\left(\left\| \Gamma_n-\Gamma\right\| _{\infty}\right). \label{tango}%
\end{equation}
The conclusion follows from the inequalities (\ref{noel})-(\ref{ppl})-(\ref{tango}), the inequality $\left\|
\Gamma_n-\Gamma\right\| _{\infty}\leq\left\| \Gamma_n-\Gamma\right\| _{2}$ and the asymptotic properties of the
sequence $(\Gamma_n-\Gamma)_{n \geq 1}$ (Bosq \cite{Bos}, Chapter 4).\hfill $\blacksquare$
\section{Some asymptotic equivalencies}
As for now, we assume $D>2$ and let
$$S_n(\bx)=\sum_{i=1}^n K\left(  \frac{\left\| \Pi_{D}\left(  \bx-\bX_{i}\right)  \right\|}{
h_n}\right) \quad \mbox{and}\quad \hat S_n(\bx)=\sum_{i=1}^n K\left(  \frac{\left\| \hat \Pi_{D}\left(  \bx-\bX_{i}\right)  \right\|}{
h_n}\right).$$
We note that $S_n(\bx)$ is a sum of independent and identically distributed random variables, whereas the terms in $\hat S_n(\bx)$ have the same distribution but are {\it not} independent. In light of the results of Section 2, our goal in this section will be to analyse the asymptotic proximity between $S_n(\bx)$ and $\hat S_n(\bx)$ under general conditions on $K$ and the sequence $(h_n)_{n \geq 1}$. Throughout, we will assume that the kernel $K$ satisfies the following set of conditions:\\

\noindent\textbf{Assumption Set $\mathbf K$}\smallskip
\begin{itemize}
\item[(${\mathbf K1}$)] $K$ is positive and bounded with compact support $[0,1]$.

\item[(${\mathbf K2}$)] $K$ is of class $\mathcal{C}^{1}$ on $[0,1]$.

\end{itemize}
These assumptions are typically satisfied by the naive kernel $K(v)=\mathbf 1_{[0,1]}(v)$. In fact, all the subsequent results also hold for kernels with an unbounded support, provided $K$ is  Lipschitz---we leave to the reader the opportunity to check the details and adapt the proofs, which turn out to be simpler in this case. For any integer $p\geq 1$, we set 
$$M_{D,p}=D\int_0^1 v^{D-1}K^{p}\left(  v\right)  \mbox{d}v$$
and, for all $\bx \in \mathcal F$ and  $h>0$, we let
$$F_{\bx}(h)=\mathbb P\left(\Pi_{D} \bX \in \mathcal{B}_{D}(\Pi_{D}\bx,h)\right),$$
where $\mathcal{B}_{D}(u,h)$ denotes the closed Euclidean ball of dimension $D$ centered at $u$ and of radius $h$. In the subsequent developments, to lighten notation a bit, and since no confusion is possible, we will write $F(h)$ instead of $F_{\bx}(h)$. Observe that $F(h)$ is positive for $\mu$-almost all $\bx \in \mathcal F$, where $\mu$ is the distribution of $\bX$. Besides, by decreasing monotonicity, since $\bX$ is nonatomic, we have
$$\lim_{h \downarrow 0} F(h)=0.$$
When the projected random variable $\Pi_D\bX$ has a density $f$ with respect to the Lebesgue measure $\lambda$ on $\mathbb  R^D$, then $F(h)\sim \gamma_{D} f(\bx)h^D$ as $h \to 0$, where $\gamma_{D}$ is a positive constant, for $\lambda$-almost all $\bx$ (see for instance Wheeden and Zygmund \cite{Wheeden}). Thus, in this case, the function $F$ is regularly varying with index $D$. We generalise this property below.\\

\noindent\textbf{Assumption Set $\mathbf R$}\smallskip
\begin{itemize}
\item[(${\mathbf R1}$)] $F$ is regularly varying at $0$ with index $D$.
\end{itemize}

Assumption $\mathbf{R}1$ means that, for any $u>0$,%
\[
\lim_{s\rightarrow0 ^+}\frac{F\left(  su\right)  }{F\left(  s\right)  }=u^{D}.%
\]
The index of regular variation was fixed to $D$ in order to alleviate the
notation, but the reader should note that our results hold for any positive index, with different constants however. In fact this index is
directly connected with the support of the distribution of $\bX$. To see this,  observe that by fixing the index to $D$ we implicitly
assume that $\Pi_{D}\bX$ fills the whole space of dimension $D$. However, elementary
calculations show that most distributions in $\mathbb{R}^{D}$,
when concentrated on a subspace of smaller dimension $D'<D$, will match assumption
$\mathbf{R}1$ with $D'$ instead of $D$. Moreover, representation theorems for
regularly varying functions (see Bingham et al. \cite{BGT}) show that, under $\mathbf R1$, $F$ may be
rewritten as $F(  u)  =u^{D}L\left(  u\right)  $, where the function $L$ is slowly
varying at $0$, that is $\lim_{s\rightarrow 0^+}L\left(  su\right)  /L\left(
s\right)  =1$. This enables to consider functions $F$ with non-polynomial behaviour such as, for instance, $F(  u)  \sim Cu^{D}|\ln u|$ as ${u\rightarrow0^{+}}$. Observe also that $F(u)$ is negligible with respect to $u^2$ as soon as $D>2$.\\

We start the analysis with two technical lemmas. Proof of Lemma \ref{L1} is deferred to Section 5, whereas Lemma  \ref{P1} is an immediate consequence of Lemma \ref{L1} and Bennett's inequality. Its proof is therefore omitted.
\begin{lem}
\label{L1} Assume that Assumption Sets $\mathbf K$ and $\mathbf R$ are satisfied. Then, for $\mu$-almost all $\bx$, if $h_n\downarrow 0$,
$$
\mathbb{E}K\left(  \frac{\left\| \Pi_{D}\left(  \bx-\bX\right)  \right\|}{
h_n}\right)   \sim M_{D,1}F\left(  h_n\right) 
$$
and
$$
\mathbb{E}K^2\left(  \frac{\left\| \Pi_{D}\left(  \bx-\bX\right)  \right\|}{h_n}\right)  \sim M_{D,2}F\left(  h_n\right)\quad \mbox{as } n \to \infty.
$$
\end{lem}
\begin{lem}
\label{P1}
Assume that Assumption Sets $\mathbf K$ and $\mathbf R$ are satisfied. Then, for $\mu$-almost all $\bx$, if $h_n\downarrow 0$ and $nF(h_n)/\ln n \to \infty$,
$$S_{n}(\bx)  \sim M_{D,1} nF\left(  h_n\right)  \quad \mbox{a.s.}$$
and
$$ \mathbb{E}S_{n}^{2}(\bx)  \sim\left[  M_{D,1}nF\left(  h_n\right)  \right]  ^{2} \quad \mbox{as } n\to \infty.$$
\end{lem}
The following proposition is the cornerstone of this section. It asserts that, asymptotically, the partial sums $S_n(\bx)$ and $\hat S_n(\bx)$ behave similarly.
\begin{pro}
\label{P2} 
Assume that Assumption Sets $\mathbf K$ and $\mathbf R$ are satisfied and that $\bX$ has bounded support. Then, for $\mu$-almost all $\bx$, if $h_n\downarrow 0$ and $nF(h_n)/\ln n \to \infty$,
\[
\hat{S}_{n}(\bx)\sim S_{n}(\bx)\quad \mbox{a.s.}
\]
and
$$\mathbb E \hat S_n^2(\bx) \sim \mathbb E S_n^2(\bx)\quad \mbox{as }n \to\infty.$$
\end{pro}
{\bf Proof of Proposition \ref{P2}}\quad To simplify notation a bit, we let, for $i=1, \hdots,n$,  $V_{i}=\|\Pi_D(\bx- \bX_{i})\|$ and $\hat V_{i}=\|\hat \Pi_D(\bx- \bX_{i})\|$. Let the events $\mathcal E_i$ and $\hat {\mathcal E_i}$ be defined by
$$
\mathcal{E}_{i}  =\left\{  V_i\leq h_n\right\} \quad \mbox{and}\quad \hat{\mathcal{E}}_{i}  =\left\{  \hat V_i\leq h_n\right\} .$$
Clearly,
\begin{align*}
&\hat S_n(\bx)-S_n(\bx)\\
&  \quad =\sum_{i=1}^{n}K\left(  \hat V_i  /{h_n}\right)  -\sum_{i=1}%
^{n}K\left(  V_i  /{h_n}\right) \\
&  \quad =\sum_{i=1}^{n}\left[  K\left( \hat  V_{i}/{h_n}\right)  -K\left( {V}%
_{i}/{h_n}\right)  \right]  \mathbf{1}_{\hat{\mathcal{ E}}_{i}\cap{\mathcal{E}}_{i}}+\sum_{i=1}^{n}K\left(  \hat V_{i}/{h_n}\right)  \mathbf{1}%
_{\hat{\mathcal{E}}_{i}\cap {{\mathcal{E}}_{i}^c}}\\
& \qquad-\sum_{i=1}^{n}K\left(  {V}_{i}/{h_n}\right)  \mathbf{1}_{{\hat{\mathcal{E}}%
_{i}^c}\cap{\mathcal{E}}_{i}}.
\end{align*}
Therefore
\begin{align}
& \left\vert \hat S_n(\bx)-S_n(\bx) \right\vert  \notag \\
& \quad \leq 
\sum_{i=1}^{n}\left\vert K\left( \hat V_{i}/{h_n}\right) -K\left( {V}%
_{i}/{h_n}\right) \right\vert \mathbf{1}_{\hat{\mathcal{E}}_{i}\cap {%
\mathcal{E}}_{i}} +C\sum_{i=1}^{n}\left( \mathbf{1}_{%
\mathcal{E}_{i}\cap {\hat{\mathcal{E}}_{i}^c}}+\mathbf{1}_{{\mathcal{E}_{i}^c}\cap \hat{\mathcal{E}}_{i}}\right)
\notag \\
& \quad \leq C \left[ \frac{\left\| \hat{\Pi }_{D}-\Pi _{D}\right\| _{\infty }}{h_n}%
\sum_{i=1}^{n}\left\| \bx -\bX_{i}\right\| \mathbf{1}_{\mathcal{E}%
_{i}}+ \sum_{i=1}^{n}\left( \mathbf{1}_{\mathcal{E}_{i}\cap {\hat{\mathcal{E}}_{i}^c}}+\mathbf{1}_{{\mathcal{E}_{i}^c}\cap \hat{\mathcal{E}}_{i}}\right)\right].\label{decomp-somme}
\end{align}%
Consequently, by Lemma \ref{P1}, the result will be
proved if we show that %
$$
\frac{\left\| \hat{\Pi }_{D}-\Pi _{D}\right\| _{\infty }}{%
n{h_n}F\left( h_n\right) }\sum_{i=1}^{n}\left\| \bx-\bX_{i}\right\| %
\mathbf{1}_{\mathcal{E}_{i}} \to 0\quad \mbox{a.s.} 
$$
and
$$
\frac{1}{nF\left( h_n\right) }\sum_{i=1}^{n}\left(\mathbf{1}_{\mathcal{E}%
_{i}\cap \hat{\mathcal{E}}_{i}^c}+\mathbf{1}_{{%
\mathcal{E}_{i}^c}\cap \hat{\mathcal{E}}_{i}}\right)  \to 0\quad
\mbox{a.s.}\quad \mbox{as } n\to \infty.
$$
The first limit is proved in technical Lemma \ref{L2} and the second one in technical Lemma \ref{L3}.\\

We proceed now to prove the second statement of the proposition. We have to show that
\begin{equation*}
\frac{\mathbb{E} \hat S^2_n(\bx)}{\mathbb{E} S_n^{2}(\bx)}\to 1\quad\mbox{as } n\to \infty.
\end{equation*}%
Using the decomposition
$$\frac{\mathbb{E}U^{2}}{\mathbb{E}V^{2}}=1+\frac{\mathbb{E}\left[
U-V\right] ^{2}}{\mathbb{E}V^{2}}+2\frac{\mathbb{E}\left[V\left( U-V\right) \right]}{%
\mathbb{E}V^{2}},$$
and the bound
$$\frac{\left|\mathbb{E}\left[V\left( U-V\right) \right]\right|}{\mathbb{E}V^{2}} \leq \sqrt{\frac{\mathbb{E}\left[
U-V\right] ^{2}}{\mathbb{E}V^{2}}},$$
it will be enough to prove that
\begin{equation*}
\frac{\mathbb{E}\left[ \hat S_n(\bx)-S_n(\bx) \right] ^{2}}{\mathbb{E}S^2_n(\bx)}\to 0,
\end{equation*}%
which in turn comes down to prove that
\begin{equation*}
\frac{\mathbb{E}\left[ \hat S_n(\bx)-S_n(\bx) \right] ^{2}}{n^{2}F^{2}\left( h_n\right) }\to 0,
\end{equation*}%
since $\mathbb{E}S_{n}^{2}(\bx)\sim \left[ M_{D,1}nF\left( h_n\right) \right] ^{2}$ by Lemma \ref{P1}.\\

Starting from inequality (\ref{decomp-somme}), we obtain
\begin{align}
& \left[ \hat S_n(\bx)-S_n(\bx) \right] ^{2} \notag\\
& \quad \leq C \left [\frac{\left\| \hat{\Pi }_{D}-\Pi _{D}\right\| _{\infty }}{%
h_n}\left( \sum_{i=1}^{n}\left\| \bx-\bX_{i}\right\| \mathbf{1}_{%
\mathcal{E}_{i}}\right)\right] ^{2}\notag\\
& \qquad +C\left[ \sum_{i=1}^{n}\left( \mathbf{1}_{\mathcal{E}_{i}\cap 
{\hat{\mathcal{E}}_{i}^c}}+\mathbf{1}_{{\mathcal{E}%
_{i}^c}\cap \hat{\mathcal{E}}_{i}}\right) \right] ^{2}\label{vacances}.
\end{align}%
Consequently, the result will be proved if we show that
$$\mathbb E \left[\frac{\left\| \hat{\Pi }_{D}-\Pi _{D}\right\| _{\infty }}{nh_nF(h_n)}\left(\sum_{i=1}^{n}\left\| \bx-\bX_{i}\right\| \mathbf{1}_{%
\mathcal{E}_{i}}\right)\right]^2\to 0.$$
and
$$\mathbb E\left[\frac{1}{nF(h_n)} \sum_{i=1}^{n}\left( \mathbf{1}_{\mathcal{E}_{i}\cap 
{\hat{\mathcal{E}}_{i}^c}}+\mathbf{1}_{{\mathcal{E}%
_{i}^c}\cap \hat{\mathcal{E}}_{i}}\right) \right] ^{2} \to 0\quad \mbox{as }n\to\infty.$$
The first limit is established in technical Lemma \ref{L4} and the second one in technical Lemma \ref{L5}.\hfill $\blacksquare$\\

The consequences of Proposition \ref{P2} in terms of kernel regression estimation will be thoroughly explored in Section 4. However, it has already important repercussions in density estimation, which are briefly sketched here and may serve as references for further studies. Suppose that the projected random variable $\Pi_D\bX$ has a density $f$ with respect to the Lebesgue measure $\lambda$ on $\mathbb  R^D$. In this case, the PCA-kernel density estimate of $f$---based on the sample $(\hat \Pi_D \bX_1, \hdots, \hat \Pi_D \bX_n)$--- reads
$$\hat f_n(\bx)=\frac{\hat{S}_{n}(\bx)}{nh_n^D}$$
and the associated pseudo-estimate---based on $(\Pi_D \bX_1, \hdots, \Pi_D \bX_n)$--- takes the form
$$f_n(\bx)=\frac{{S}_{n}(\bx)}{nh_n^D}.$$
An easy adaptation of the proof of Corollary \ref{C1} in Section 4 shows that, under the conditions of Proposition \ref{P2},
 $$\mathbb{E}%
\left[ \hat{f}_n(\bx) -{f_n}(\bx) \right] ^{2}=\mathcal O\left( \frac{\log \left( nh_n^{2}\right) }{nh_n^{2}}\right).$$

To illustrate the importance of this result, suppose for example that the target density $f$ belongs to the class $\mathcal G_p$ of $p$-times continuously differentiable functions. In this context (Stone \cite{Stone1, Stone2}), the optimal rate of convergence over $\mathcal G_p$ is $%
n^{-2p/\left( 2p+D\right) }$ and the kernel density estimate with a bandwidth $h_n^{\ast }\asymp n^{-1/\left( 2p+D\right)
}$ achieves this minimax rate. Thus, letting $\hat{f}_n^{\ast }$ (respectively ${f}_n^{\ast})$ be the PCA-kernel density estimate (respectively pseudo-density estimate) based on this optimal bandwidth, we are led to%
\[
\frac{\mathbb{E}\left[ \hat{f}_n^{\ast }(\bx) -{f}_n%
^{\ast }( \bx) \right] ^{2}}{n^{-2p/\left( 2p+D\right)}}\rightarrow 0
\]%
as soon as $D>2$. Thus, the $L_2$-rate of convergence of $\hat{f}_n^{\star}$ towards $f_n^{\star}$ is negligible with respect to the $L_2$-rate of convergence of $f_n^{\star}$ towards $f$. In consequence, replacing $\hat f_n^{\star}$ by $f_n^{\star}$ has no effect on the asymptotic rate. The same ideas may be transposed without further effort to asymptotic normality and other error criteria.
\section{Regression analysis}
As framed in the introduction, we study in this final section the PCA-kernel regression procedure, which was our initial motivation. Recall that, in this context, we observe a set $\{(\bX_1,Y_1), \hdots, (\bX_n,Y_n)\}$ of independent $\mathcal F\times \mathbb R$-valued random variables with the same distribution as a generic pair $(\bX, Y)$, where $\bX$ is nonatomic centered, and $Y$ satisfies $\mathbb E |Y|<\infty$. The goal is to estimate the regression function $r^{D}(\bx)=\mathbb E[Y|\Pi_D\bX=\Pi_D\bx]$ via the PCA-kernel estimate, which takes the form
$$\hat r_n^D(\bx)=\frac{\sum_{i=1}^n Y_iK\left(\frac{\|\hat \Pi_D(\bx-\bX_i)\|}{h_n} \right)}{\sum_{i=1}^nK\left(\frac{\|\hat \Pi_D(\bx-\bX_i)\|}{h_n} \right)}.$$
This estimate is mathematically intractable and we plan to prove that we can substitute without damage to $\hat r_n^D$ the pseudo-estimate
$$r_n^D(\bx)=\frac{\sum_{i=1}^n Y_iK\left(\frac{\|\Pi_D(\bx-\bX_i)\|}{h_n} \right)}{\sum_{i=1}^nK\left(\frac{\|\Pi_D(\bx-\bX_i)\|}{h_n} \right)}.$$
To this aim, observe first that, with the notation of Section 3,
$$\hat r_n^D(\bx)=\frac{\hat Z_n(\bx)}{\hat S_n(\bx)}$$
and
$$r_n^D(\bx)=\frac{Z_n(\bx)}{S_n(\bx)},$$
where, for all $n\geq 1$, 
$$
\hat{Z}_{n}(\bx)=\sum_{i=1}^{n}Y_{i}K\left( \frac{\left\| \hat{\Pi }_{D}\left(\bx- \bX_{i}\right) \right\|}{h_n}\right)
$$
and
$$Z_{n}(\bx)=\sum_{i=1}^{n}Y_{i}K\left( \frac{\left\| \Pi _{D}\left(\bx-\bX_{i}\right)
\right\|}{h_n}\right).$$
\begin{pro}
\label{T2}
Assume that Assumption Sets $\mathbf K$ and $\mathbf R$ are satisfied, that $\bX$ has bounded support and $Y$ is bounded. Assume also that $r^{D}(\bx)\neq0$ and $r^{D}$ is Lipschitz in a neighborhood of $\bx$. Then, for $\mu$-almost all $\bx$, if $h_n\downarrow 0$ and $nF(h_n)/\ln n \to \infty$,

$$\hat{Z}_{n}(\bx)\sim Z_{n}(\bx)\quad \mbox{a.s.}$$
and 
$$\mathbb{E}\hat{Z}_{n}^{2}(\bx)\sim \mathbb{E}Z_{n}^{2}(\bx)\quad \mbox{as } n\to \infty.$$
\end{pro}
\noindent{\bf Proof of Proposition \ref{T2}}\quad  Using the Lipschitz property of $r^D$,  we easily obtain by following the lines of Lemma \ref{L1} that, at $\mu$-almost all $\bx$, 
$$\mathbb{E}\left[YK\left(\frac{\left\| \Pi _{D}\left( \bx-\bX\right) \right\| }{h_n}\right)\right] \sim r^{D}( \bx) F(h_n)$$
and
$$\mathbb{E}\left[Y^{2}K^{2}\left(\frac{\left\| \Pi _{D}\left( \bx-\bX\right) \right\| }{h_n}\right)\right] \sim C F\left(h_n\right) \quad \mbox{as } n \to \infty.$$
Moreover, for $\mu$-almost all $\bx$,  
$$Z_{n}(\bx) \sim n\mathbb{E}\left[YK\left(
\frac{\left\| \Pi _{D}\left( \bx-\bX\right) \right\| }{h_n}\right)\right]\quad \mbox{a.s.}$$
and
$$
\mathbb E Z^2_{n}(\bx) \sim \left [n\mathbb{E}\left[YK\left(
\frac{\left\| \Pi _{D}\left( \bx-\bX\right) \right\| }{h_n}\right)\right]\right]^2\quad  \mbox{as } n \to \infty.$$
The first equivalence is a consequence of Bennett's inequality and the fact that, for all large enough $n$ and $i =1,2$,
$$\mathbb{E}\left[\left|Y\right|^{i}K\left(
\frac{\left\| \Pi _{D}\left( \bx-\bX\right) \right\| }{h_n}\right)\right] \geq C F\left( h_n\right),$$
which itself follows from the requirement $r^D(\bx)>0$.\\

Finally, since $Y$ is bounded, an inspection of the proof of Proposition \ref{P2} reveals that displays (\ref{decomp-somme}) and (\ref{vacances}) may be verbatim repeated with $S$ replaced by $Z$.\hfill $\blacksquare$
\begin{cor}
\label{C1} Under the assumptions of Proposition \ref{T2}, for $\mu$-almost all $\bx$, the estimate $\hat{r}_n^D$ and the pseudo-estimate $r^{D}_n$ satisfy
\begin{equation*}
{\hat{r}_n^D\left( \bx\right) } \sim {r^{ D }_n\left( \bx\right) }\quad \mbox{a.s. as } n \to \infty.
\end{equation*}%
Moreover,%
\begin{equation*}
\mathbb{E}\left[ \hat{r}_n^D\left( \bx\right) -r^D_n\left(\bx\right) \right] ^{2}=\mathcal {O}\left( \frac{\log\left(nh_n^{2}\right)}{nh_n^{2}}\right).
\end{equation*}%
\end{cor}
\noindent {\bf Proof of Corollary \ref{C1}}\quad We start with the decomposition
\begin{equation*}
\hat{r}_n^D\left( \bx\right) -r^{D }_n\left( \bx\right) =\frac{\hat{%
Z}_{n}(\bx)}{\hat{S}_{n}(\bx)}\left( 1-\frac{\hat{S}_{n}(\bx)}{S_{n}(\bx)}\right) +\frac{%
1}{S_{n}(\bx)}\left( \hat{Z}_{n}(\bx)-Z_{n}(\bx)\right),
\end{equation*}%
which comes down to
\begin{equation*}
\frac{\hat{r}_n^D\left( \bx\right) -r^{D }_n\left( \bx\right)}{\hat{r}_n^D\left( \bx\right)}=\left( 1-\frac{\hat{S}_{n}(\bx)}{S_{n}(\bx)}\right)+\frac{\hat{S}_{n}(\bx)}{S_{n}(\bx)}\left( 1-\frac{Z_{n}(\bx)}{\hat{Z}_{n}(\bx))}\right).
\end{equation*}
The first part of the corollary is then an immediate consequence of Proposition \ref{T1} and Proposition \ref{T2}.\\

We turn to the second part. Note that $\hat{Z}_{n}(\bx)/\hat{S}%
_{n}(\bx)$ is bounded whenever $Y$ is bounded. In consequence, we just need to provide upper bounds for the terms $\mathbb{E}\left[ 1-\hat{S}_{n}(\bx)/S_{n}%
(\bx)\right] ^{2}$ and $\mathbb{E}\left[ \left( \hat{Z}_{n}(\bx)-Z_{n}(\bx)\right)
/S_{n}(\bx)\right] ^{2}.$  Besides, classical arguments show that the latter two expectations may be replaced by $\mathbb{E}\left[
S_{n}(\bx)-\hat{S}_{n}(\bx)\right] ^{2}/\mathbb{E}S^2_{n}(\bx)$ and $%
\mathbb{E}\left[ \hat{Z}_{n}(\bx)-Z_{n}(\bx)\right] ^{2}/\mathbb{E}S_{n}^{2}(\bx)$, respectively. It turns out that the analysis of each of these terms is similar, and we will therefore focus on the first one only. Given the result of Proposition \ref{P1}, this comes down to analyse $\mathbb{E}\left[ S_{n}(\bx)-\hat{S}%
_{n}(\bx)\right] ^{2}/\left[ M_{D,1}nF\left( h_n\right) \right] ^{2}$ and to refine the bound. \\

By inequality (\ref{vacances}), we have %
\begin{align}
\frac{\mathbb E\left[ S_{n}(\bx)-\hat{S}_{n}(\bx)\right]^2}{n^2F^2(h_n)}
& \quad \leq C\mathbb E\left[\frac{\left\| \hat{\Pi }_{D}-\Pi _{D}\right\| _{\infty }}{nh_nF(h_n)}\left( \sum_{i=1}^{n}\left\| \bx-\bX_{i}\right\| \mathbf{1}_{%
\mathcal{E}_{i}}\right)\right] ^{2}\notag\\
& \qquad +C\left[ \frac{1}{nF(h_n)}\sum_{i=1}^{n}\left( \mathbf{1}_{\mathcal{E}_{i}\cap 
{\hat{\mathcal{E}}_{i}^c}}+\mathbf{1}_{{\mathcal{E}%
_{i}^c}\cap \hat{\mathcal{E}}_{i}}\right) \right] ^{2}\label{vacances2}.
\end{align}%
With respect to the first term, Lemma \ref{L4} asserts that
\begin{equation*}
\mathbb{E}\left[ \frac{\left\| \hat{\Pi }_{D}-\Pi _{D}\right\|
_{\infty }}{nh_nF(h_n)}\left(\sum_{i=1}^{n}\left\|
\bx-\bX_{i}\right\| \mathbf{1}_{\mathcal{E}_{i}}\right)\right] ^{2}=\mathcal {O}\left( 
\frac{1}{nh_n^{2}}\right).
\end{equation*}%
The second term in inequality (\ref{vacances2}) is of the order $\mathcal O (\log (nh_n^2)/(nh_n^2))$, as proved in technical Lemma \ref{L6}. This completes the proof.\hfill $\blacksquare$\\

To illustrate the usefulness of Corollary \ref{C1}, suppose that the regression function $r^D$ belongs to the class $\mathcal G_p$ of $p$-times continuously differentiable functions.
In this framework, it is well-known (Stone \cite{Stone1, Stone2}) that the optimal rate of convergence on the class $\mathcal G_p$ is $%
n^{-2p/\left( 2p+D\right) }$ and that the kernel estimate with a bandwidth $h_n^{\ast }\asymp n^{-1/\left( 2p+D\right)
}$ achieves this minimax rate. Plugging this optimal $h_n^{\ast }$ into the rate of Corollary \ref{C1}, we obtain%
\begin{equation*}
\mathbb{E}\left[ \hat{r}_n^D\left( \bx\right) -r_n^D\left(
\bx\right) \right] ^{2}=\mathcal{O}\left( n^{-\alpha }\right),
\end{equation*}
with $\alpha =\left( 2p+D-2\right) /\left( 2p+D\right)$. This rate is strictly faster than the minimax rate $n^{-2p/\left( 2p+D\right) }$ provided $%
\alpha >2p/\left( 2p+D\right) $ or, equivalently, when $D>2$. In this case, 
\begin{equation*}
\lim_{n\to \infty }\frac{\mathbb{E}\left[ \hat{r}_n^D\left(
\bx\right) -r^{D}_n\left( \bx\right) \right] ^{2}}{n^{-2p/\left( 2p+D\right) }}=0,
\end{equation*}
and Corollary \ref{C1} claims in fact that the rate of convergence of $\hat{r}_n^D $ towards $r^{D}_n$ is negligible with respect to the rate of convergence of $r_n^D$ towards $r^D$. In conclusion, even if $\hat{r}_n^D$ is the only
possible and feasible estimate, carrying out its asymptotics from the pseudo-estimate $r^{D
}_n$ is permitted.
\section{Proofs}
\subsection{Proof of Lemma \ref{L0}}
Observe first, since $\hat{\lambda}_{1}\geq\hdots\geq\hat{\lambda}_{D}$, that
$$\mathcal{A}_{n}=\left\{ \hat{\lambda}_1\leq \lambda_1+1/2, \hat{\lambda}_{D}\geq \lambda_{D}-\delta
_{D},\mbox{and } \hat{\lambda}_{D+1}<\lambda_{D}-\delta_{D}\right\}.$$
Therefore
\begin{align*}
& \mathbb{P}\left({\mathcal{A}}^c_{n}\right) \\
  &\quad  \leq\mathbb{P}\left(
\hat{\lambda}_{1}-\lambda_{1}>1/2\right)+\mathbb{P}\left(
\hat{\lambda}_{D}-\lambda_{D}<-\delta_{D}\right)  +\mathbb{P}\left(
\hat{\lambda}_{D+1}-\lambda_{D+1}\geq\delta_{D}\right) \\
& \quad  \leq\mathbb{P}\left(  \left\vert
\hat{\lambda}_{1}-\lambda_{1}\right\vert \geq1/2\right)+\mathbb{P}\left(  \left\vert \hat{\lambda}_{D}-\lambda
_{D}\right\vert \geq\delta_{D}\right)  +\mathbb{P}\left(  \left\vert
\hat{\lambda}_{D+1}-\lambda_{D+1}\right\vert \geq\delta_{D}\right).
\end{align*}
The inequality %
\[
\sup_{i \geq 1}\left\vert \hat{\lambda}_{i}-\lambda_{i}\right\vert \leq\left\|
\Gamma_{n}-\Gamma\right\| _{2}%
\]
shifts the problem from $\vert \hat{\lambda
}_{i}-\lambda_{i}\vert $ to $\| \Gamma_{n}-\Gamma\|
_{2}.$ An application of a standard theorem for Hilbert-valued random variables
(see for instance Bosq \cite{Bos}) leads to
\[
\mathbb{P}\left(  \left\| \Gamma_{n}-\Gamma\right\| _{2}\geq \varepsilon
\right)  =\mathcal O\left(\exp\left(  -c_{1}\frac{n\varepsilon^{2}}{c_{2}%
+c_{3}\varepsilon}\right)\right),
\]
for three positive constants $c_1$, $c_2$ and $c_3$. Consequently, for fixed $D$,
\[
\mathbb{P}\left(  \mathcal{A}_{n}^{c}\right)  =\mathcal O\left(\exp\left(
-nC\varepsilon^2\right)\right),
\]
where $C$ is a positive constant depending on $D$.
\subsection{Proof of Lemma \ref{L1}}
The proof will be based on successive applications of Fubini's theorem. Denoting by $\mu_{D,\bx,h_n}$ the probability measure associated with the random variable
$\left\| \Pi_{D}\left(  \bx-\bX\right)  \right\| /h_n$, we may write
\[
\mathbb{E}K\left(  \frac{\left\| \Pi_{D}\left(  \bx-\bX\right)  \right\|}{
h_n}\right) =\int_{0}^{1}K\left(  v\right)  \mu_{D,\bx,h_n}\left( \mbox{d}v\right).
\]
Thus
\begin{align*}
\mathbb{E}K\left(  \frac{\left\| \Pi_{D}\left(  \bx-\bX\right)  \right\|}{
h_n}\right)  &  =\int_{0}^{1}\left[  K\left(  1\right)  -\int_{v}^{1}K^{\prime
}\left(  s\right) \mbox{d}s\right]  \mu_{D,\bx,h_n}\left(  \mbox{d}v\right)  \\
&  =K\left(  1\right)  F\left(  h_n\right)  -\int_0^1K^{\prime}\left(
s\right) \int_{[0\leq v\leq s] }   \mu_{D,\bx,h_n}
\left(\mbox{d}v\right)\mbox{d}s  \\
&  =K\left(  1\right)  F\left(  h_n\right)  -\int_{0}^{1}F\left(  h_ns\right)
K^{\prime}\left(  s\right)  \mbox{d}s\\
&  =F\left(  h_n\right)  \left[ K\left(  1\right)  -\int_{0}^{1}\frac{F\left(
h_ns\right)  }{F\left(  h_n\right)  }K^{\prime}\left(  s\right)  \mbox{d}s\right].
\end{align*}
Using the fact that $F$ is increasing regularly varying of order $D$, an application of Lebegue's dominated
convergence theorem yields%
$$
\mathbb{E}K\left(  \frac{\left\| \Pi_{D}\left(  \bx-\bX\right)  \right\|}{
h_n}\right) \sim F\left(  h_n\right)  \left[  K\left(  1\right)  -\int_{0}%
^{1}s^{D}K^{\prime}\left(  s\right)  \mbox{d}s\right] 
$$
i.e.,
$$
\mathbb{E}K\left(  \frac{\left\| \Pi_{D}\left(  \bx-\bX\right)  \right\|}{
h_n}\right)\sim F\left(  h_n\right)  D\int
_{0}^{1}s^{D-1}K\left(  s\right) \mbox{d}s\quad \mbox{as } n\to \infty.
$$
This shows the first statement of the lemma. Proof of the second statement is similar.
\subsection{Some technical lemmas}
In this subsection, for all $i=1, \hdots, n$, we let $V_{i}=\|\Pi_D(\bx- \bX_{i})\|$ and $\hat V_{i}=\|\hat \Pi_D(\bx- \bX_{i})\|$. The events $\mathcal E_i$ and $\hat {\mathcal E_i}$ are defined by
$$
\mathcal{E}_{i}  =\left\{  V_i\leq h_n\right\} \quad \mbox{and}\quad \hat{\mathcal{E}}_{i}  =\left\{  \hat V_i\leq h_n\right\} .$$
\begin{lem}
\label{L2}%
Assume that $\bX$ has bounded support. Then,  for $\mu$-almost all $\bx$, 
\begin{equation*}
\frac{\left\| \hat{\Pi}_{D}-\Pi_{D}\right\| _{\infty}}{nh_nF\left(
h_n\right) }\sum_{i=1}^{n}\left\| \bx-\bX_{i}\right\| \mathbf{1}_{\mathcal{E}_{i}}=\mathcal O \left( \sqrt{\frac{\log n}{nh_n^2}}\right)\quad \mbox{a.s.}
\end{equation*}
\end{lem}
{\bf Proof of Lemma \ref{L2}}\quad According to statement $(ii)$ of Theorem \ref{T1}, 
$$\frac{\left\| \hat{\Pi}_{D}-\Pi_{D}\right\| _{\infty}}{nh_nF(h_n)}=\mathcal {O}\left( \sqrt{\frac{\log n}{%
n^3h_n^{2}F^2(h_n)}}\right) \quad \mbox{a.s.} 
$$
Moreover, since $\bX$ has bounded support, there exists a positive constant $M$ such that, for $\mu$-almost all $\bx$, 
$$
\sum_{i=1}^{n}\left\| \bx-\bX_{i}\right\| \mathbf{1}_{\mathcal{E}_{i}} \leq M\sum_{i=1}^{n}\mathbf{1}_{\mathcal{E}_{i}}\quad \mbox{a.s.}$$
Clearly, $\sum_{i=1}^{n}\mathbf{1}_{\mathcal{E}_{i}}$ has a Binomial distribution with parameters $n$ and $F\left(h_n\right)$ and consequently, by Bennett's inequality,
$$\sum_{i=1}^{n}\left\| \bx-\bX_{i}\right\| \mathbf{1}_{\mathcal{E}_{i}}=\mathcal O\left(nF\left(h_n\right)\right) \quad \mbox{a.s.}$$
This completes the proof of the lemma.\hfill $\blacksquare$
\begin{lem}
\label{L3}
Assume that Assumption Set $\mathbf R$ is satisfied and $\bX$ has bounded support. Then, if $h_n\downarrow 0$ and $nh_n^2/\log n \to \infty$,
\begin{equation*}
\label{eq-ind}
\frac{1}{nF\left( h_n\right) }\sum_{i=1}^{n}\left(\mathbf{1}_{\mathcal{E}%
_{i}\cap \hat{\mathcal{E}}_{i}^c}+\mathbf{1}_{{%
\mathcal{E}_{i}^c}\cap \hat{\mathcal{E}}_{i}}\right)  \to 0\quad
\mbox{a.s.}\quad \mbox{as } n\to \infty.
\end{equation*}
\end{lem}

{\bf Proof of Lemma \ref{L3}}\quad Define $\kappa_{n}=C_{\kappa} \sqrt{\frac{\log n}{nh_n^2}}$ and $\eta_n= \kappa_n h_n$, where $C_{\kappa}$ is a constant which will be chosen later. Observe that%
\begin{equation*}
\mathcal{E}_{i}\cap {\hat{\mathcal{E}}_{i}^c}= \left[\left\{ h_n-\eta_n
<V_{i}\leq h_n\right\} \cap {\hat{\mathcal{E}}_{i}^c}\right]\cup\left[\left\{ V_{i}\leq
h_n-\eta_n \right\} \cap {\hat{\mathcal{E}}_{i}^c}\right].
\end{equation*}%
Similarly
\begin{equation*}
{\hat{\mathcal{E}}}_{i}^c=\left\{ \hat{V}_{i}>h_n\right\}
=\left\{ \hat{V}_{i}-V_{i}>h_n-V_{i}\right\}.
\end{equation*}%
Consequently, we may write
\begin{align}
\sum_{i=1}^{n}\mathbf{1}_{\mathcal{E}_{i}\cap {\hat{\mathcal{%
E}}_{i}^c}}& \leq \sum_{i=1}^{n}\mathbf{1}_{\left\{ h_n-\eta_n <V_{i}\leq
h_n\right\} }+\sum_{i=1}^{n}\mathbf{1}_{\left\{ V_{i}\leq h_n-\eta_n \right\}
\cap {\hat{\mathcal{E}}_{i}^c}}  \notag \\
& =\sum_{i=1}^{n}\mathbf{1}_{\left\{ h_n-\eta_n <V_{i}\leq h_n\right\}
}+\sum_{i=1}^{n}\mathbf{1}_{\left\{ V_{i}\leq
h_n-\eta_n, \hat{V}_{i}-V_{i}>h_n-V_i \right\} }  \notag \\
& \leq \sum_{i=1}^{n}\mathbf{1}_{\left\{ h_n-\eta_n <V_{i}\leq h_n\right\}
}+\sum_{i=1}^{n}\mathbf{1}_{\left\{ \hat{V}_{i}-V_{i}>\eta_n \right\} }\notag\\
& \leq \sum_{i=1}^{n}\mathbf{1}_{\left\{ h_n-\eta_n <V_{i}\leq h_n\right\}
}+\sum_{i=1}^{n}\mathbf{1}_{\left\{ \left|\hat{V}_{i}-V_{i}\right|>\eta_n \right\} }.\label{majo-indic}
\end{align}%
By Bennett's inequality, we have %
\begin{equation*}
\sum_{i=1}^{n}\mathbf{1}_{\left\{ h_n-\eta_n <V_{i}\leq h_n\right\} }\sim n%
\mathbb{P}\left( h_n-\eta_n <V_1 \leq h_n\right)\quad \mbox{a.s. as } n \to \infty,
\end{equation*}%
whence%
\begin{equation}
\label{train}
\frac{1}{nF\left( h_n\right) }\sum_{i=1}^{n}\mathbf{1}_{\left\{ h_n-\eta_n
<V_{i}\leq h_n\right\} }\sim \frac{F\left( h_n\right) -F\left( h_n-\eta_n \right) }{%
F\left( h_n\right) }\quad \mbox{a.s.}
\end{equation}
But, using the fact that $F$ is regularly varying with index $D$, we may write
\begin{equation}
\label{gg}
1-\frac{F\left( h_n-\eta_n \right) }{F\left( h_n\right)}\sim 1-\left(1-\kappa_n\right)^{D}\sim  D \kappa_n \to 0 \quad \mbox{as }n \to \infty.
\end{equation}
(Note that the value of the index $D$ influences constants only.)\\

Next, since $\bX$  has bounded support, at $\mu$-almost all $\bx$, %
\begin{align*}
 \sum_{i=1}^{n}\mathbf{1}_{\left\{ \left\vert \hat{V}%
_{i}-V_{i}\right\vert >\eta_n \right\} }& \leq \sum_{i=1}^{n}\mathbf{1}%
_{\left\{ \left\| \hat{\Pi }_{D}-\Pi _{D}\right\| _{\infty
}\left\| \bx-\bX_{i}\right\| >\eta_n \right\} }\\
&\leq n\mathbf{1}_{\left\{ \left\| \hat{\Pi }_{D}-\Pi _{D}\right\|
_{\infty }>\eta_n /M\right\} } \quad \mbox{a.s.}
\end{align*}%
for some positive constant $M$. By statement $(i)$ of Theorem \ref{T1}, we have
\begin{equation*}
\mathbb{P}\left( \left\| \hat{\Pi }_{D}-\Pi _{D}\right\| _{\infty
}\geq \varepsilon \right) =\mathcal {O}\left(\exp \left( -Cn\varepsilon^{2} \right)\right).
\end{equation*}%
Therefore,
\begin{equation}
\label{bb}
\frac{1}{nF(h_n)}\sum_{i=1}^{n}\mathbf{1}_{\left\{  \left|\hat{V}_{i}-V_{i}\right| >\eta_n \right\}}\to
0\quad \mbox{a.s.}
\end{equation}%
whenever
$$\sum_{n=1}^{\infty }\exp \left( -Cn\eta_n^{2} /M^2\right)
<\infty.$$
Observing that $n\eta_n^{2}= C_{\kappa}^2 \log n$,  we see that the summability condition above is fulfilled as soon as $C_{\kappa}$ is large enough.\\

Combining inequality (\ref{majo-indic}) with (\ref{train})-(\ref{gg}) and (\ref{bb}), we conclude that
$$\frac{1}{nF(h_n)}\sum_{i=1}^n\mathbf{1}_{\mathcal{E}_{i}\cap {\hat{\mathcal{%
E}}_{i}^c}} \to 0 \quad \mbox{a.s. as } n \to \infty.$$
One shows with similar arguments that 
$$\frac{1}{nF(h_n)}\sum_{i=1}^{n}\mathbf{1}_{{\mathcal{E}_{i}^c}\cap 
\hat{\mathcal{E}}_{i}} \to 0 \quad \mbox{a.s. as } n \to \infty.$$
$\hfill \blacksquare$
\begin{lem}
\label{L4}
Assume that $\bX$ has bounded support. Then,  for $\mu$-almost all $\bx$, if 
$nF\left(  h_n\right)  \to \infty$,
$$\mathbb E \left[\frac{\left\| \hat{\Pi }_{D}-\Pi _{D}\right\| _{\infty }}{nh_nF(h_n)}\left(\sum_{i=1}^{n}\left\| \bx-\bX_{i}\right\| \mathbf{1}_{%
\mathcal{E}_{i}}\right)\right]^2=\mathcal O \left ( \frac{1}{nh_n^2}\right).$$
\end{lem}
{\bf Proof of Lemma \ref{L4}}\quad  For $i=1, \hdots, n$, let $U_{i}=\left\| \bx-\bX_{i}\right\| \mathbf{1}_{\mathcal{E}%
_{i}}-\mathbb{E}\left[\left\| \bx-\bX_{i}\right\| \mathbf{1}_{\mathcal{E}%
_{i}}\right],$
and write
\begin{align*}
& \mathbb{E}\left[ \left\| \hat{\Pi }_{D}-\Pi
_{D}\right\| _{\infty }\left( \sum_{i=1}^{n}\left\|
\bx-\bX_{i}\right\| \mathbf{1}_{\mathcal{E}_{i}}\right)\right]^2 \\
& \quad \leq 2\mathbb{E}\left[ \left\| \hat{\Pi }_{D}-\Pi
_{D}\right\| _{\infty } \sum_{i=1}^{n}U_{i}\right]^2
\\
& \qquad +2n^{2} \mathbb{E}^2\left[\left\| \bx-\bX\right\| 
\mathbf{1}_{\mathcal{E}_{1}}\right] \mathbb{E}\left\| \hat{%
\Pi }_{D}-\Pi _{D}\right\| _{\infty }^{2}.
\end{align*}%
By Theorem \ref{T1} $(iii)$, 
$$\mathbb{E}\left\| \hat{\Pi }_{D}-\Pi _{D}\right\|
_{\infty }^{2}=\mathcal{O}\left( \frac{1}{n}\right), $$
and, since $\bX$ has bounded support, for $\mu$-almost all $\bx$,
 $$ \mathbb{E}^2\left[\left\|
\bx-\bX\right\| \mathbf{1}_{\mathcal{E}_{1}}\right]=\mathcal {O}\left(
F^{2}\left( h_n\right) \right).$$ 
Consequently,%
\begin{equation*}
n^2\mathbb{E}^2\left[\left\| \bx-\bX\right\| \mathbf{1}_{\mathcal{E}_{1}}\right]\mathbb{E}\left\| \hat{%
\Pi }_{D}-\Pi _{D}\right\| _{\infty }^{2}=\mathcal {O}\left(nF^{2}(h_n)\right).
\end{equation*}
One easily shows, with methods similar to the ones used to prove $(iii)$ of Theorem \ref{T1}, that
$$\mathbb{E} \left\|
\hat{\Pi}_{D}-\Pi_{D}\right\| _{\infty}^{4}=\mathcal O \left (\frac{1}{n^{2}}\right).$$
Moreover, simple computations lead to
\[
\mathbb{E}\left [ \sum_{i=1}^{n}U_{i} \right] ^{4}=\mathcal O \left(  nF\left(
h_n\right)  +n^{2}F^{2}\left(  h_n\right)  \right)  =\mathcal O\left(  n^{2}F^{2}\left(
h_n\right)  \right),
\]
when $nF\left(h_n\right)  \rightarrow \infty$. Consequently, by Cauchy-Schwarz inequality,
\begin{align*}
\mathbb{E}\left[  \left\| \hat{\Pi}_{D}-\Pi_{D}\right\| _{\infty
} \sum_{i=1}^{n}U_{i}\right]^2   &  =\mathcal O\left(  \frac
{1}{n}\left(  \mathbb{E}\left[  \sum_{i=1}^{n}U_{i}\right]  ^{4}\right)
^{1/2}\right)  \\
&  =\mathcal O\left(  F\left(  h_n\right)  \right)  .
\end{align*}

Putting all the pieces together, we obtain
$$\mathbb E \left[\frac{\left\| \hat{\Pi }_{D}-\Pi _{D}\right\| _{\infty }}{nh_nF(h_n)}\left(\sum_{i=1}^{n}\left\| \bx-\bX_{i}\right\| \mathbf{1}_{%
\mathcal{E}_{i}}\right)\right]^2=\mathcal O \left ( \frac{1}{nh_n^2}\right)\quad \mbox{as } n\to \infty.$$
\hfill $\blacksquare$
\begin{lem}
\label{L5}
Assume that Assumption Set $\mathbf R$ is satisfied and $\bX$ has bounded support. Then, if $h_n \downarrow 0$ and $nh_n^2/\log n \to \infty$, 
$$\mathbb E\left[\frac{1}{nF(h_n)} \sum_{i=1}^{n}\left( \mathbf{1}_{\mathcal{E}_{i}\cap 
{\hat{\mathcal{E}}_{i}^c}}+\mathbf{1}_{{\mathcal{E}%
_{i}^c}\cap \hat{\mathcal{E}}_{i}}\right) \right] ^{2} \to 0\quad \mbox{as } n\to \infty.$$
\end{lem}
{\bf Proof of Lemma \ref{L5}}\quad The proof is close to the derivation of Lemma \ref{L3}---almost sure convergence is replaced here by convergence in mean square. Therefore, we go quickly through it. \\

Because of $(a+b)^2\leq 2a^2+2b^2$, it is enough to prove that
$$\mathbb E\left[\frac{1}{nF(h_n)} \sum_{i=1}^{n} \mathbf{1}_{\mathcal{E}_{i}\cap 
{\hat{\mathcal{E}}_{i}^c}} \right] ^{2} \to 0$$
and
$$\mathbb E\left[ \frac{1}{nF(h_n)}\sum_{i=1}^{n} \mathbf{1}_{\mathcal{E}_{i}^c\cap 
{\hat{\mathcal{E}}_{i}}} \right] ^{2} \to 0 \quad \mbox{as }n \to \infty.$$
We will focus on the first limit only---proof of the second one is similar. With the notation of Lemma \ref{L3}, by inequality (\ref{majo-indic}),%
\begin{align}
\mathbb E \left[\frac{1}{nF(h_n)}\sum_{i=1}^{n}\mathbf{1}_{\mathcal{E}_{i}\cap {\hat{\mathcal{%
E}}_{i}^c}}\right]^2&\leq 2\mathbb E \left[\frac{1}{nF(h_n)}\sum_{i=1}^{n}\mathbf{1}_{\left\{ h_n-\eta_n <V_{i}\leq
h_n\right\} }\right]^2 \notag\\
& \quad +2\mathbb E \left[\frac{1}{nF(h_n)} \sum_{i=1}^{n}\mathbf{1}_{\left\{ \left\vert \hat{V}%
_{i}-V_{i}\right\vert >\eta_n \right\} }\right]^2,\label{hermes}
\end{align}%
where $\eta_n $ is a tuning parameter which will be fixed later.\\

The first term on the right of (\ref{hermes}) is handled exactly as in Lemma \ref{L3} and tends to zero. We just require that $\eta_n=\kappa_nh_{n}$, with $\kappa_n \to 0$.\\

With respect to the second term, write%
\begin{align*}
\mathbb{E}\left[\frac{1}{nF(h_n)} \sum_{i=1}^{n}\mathbf{1}_{\left\{ \left\vert \hat{V}%
_{i}-V_{i}\right\vert >\eta_n \right\} }\right] ^{2}& \leq \mathbb{E}\left[\frac{1}{nF(h_n)}
\sum_{i=1}^{n}\mathbf{1}_{\left\{ \left\| \hat{\Pi }_{D}-\Pi
_{D}\right\| _{\infty }\left\| \bx- \bX_{i}\right\| >\eta_n \right\}}\right] ^{2} \\
& \leq \mathbb{E}\left[\frac{1}{nF(h_n)}
\sum_{i=1}^{n}\mathbf{1}_{\left\{ \left\| \hat{\Pi }_{D}-\Pi
_{D}\right\| _{\infty } >\eta_n/M \right\}}\right] ^{2} \\
& \leq \frac{1}{F^2(h_n)}\mathbb{P}\left( \left\| \hat{\Pi }_{D}-\Pi
_{D}\right\| _{\infty }>\eta_n /M\right),
\end{align*}%
at $\mu$-almost all $\bx$ and for some positive constant $M$. Applying finally statement $(i)$ of Theorem \ref{T1}, we obtain
\begin{equation*}
\mathbb{E}\left[\frac{1}{nF(h_n)} \sum_{i=1}^{n}\mathbf{1}_{\left\{ \left\vert \hat{V}%
_{i}-V_{i}\right\vert >\eta_n \right\} }\right] ^{2}=\mathcal O\left(\frac{\exp \left( -Cn \kappa_n^{2} /M^{2} \right) }{F^{2}\left( h_n\right) }\right)
\end{equation*}%
which tends to zero whenever $\kappa_{n}=C_{\kappa} \sqrt{\frac{\log n}{nh_n^2}}$ for a sufficiently large $C_{\delta}$.\hfill $\blacksquare$
\begin{lem}
\label{L6}
Assume that Assumption Set $\mathbf R$ is satisfied and $\bX$ has bounded support. Then, if $nF(h_n) \to \infty$, 
$$\mathbb E\left[\frac{1}{nF(h_n)} \sum_{i=1}^{n}\left( \mathbf{1}_{\mathcal{E}_{i}\cap 
{\hat{\mathcal{E}}_{i}^c}}+\mathbf{1}_{{\mathcal{E}%
_{i}^c}\cap \hat{\mathcal{E}}_{i}}\right) \right] ^{2} = \mathcal O \left( \frac{\log (nh_n^2)}{nh_n^2}\right).$$
\end{lem}
{\bf Proof of Lemma \ref{L6}}\quad  We deal only with the term
$$\mathbb E \left [\frac{1}{nF(h_n)}\sum_{i=1}^{n}\mathbf{1}_{\mathcal{E}_{i}\cap {%
\hat{\mathcal{E}}_{i}}^c}\right]^2,$$
since the other one may be addressed the same way. At this
point, we have to get sharper into the bounds derived in Lemma \ref{L5}. Let $(\kappa_n)_{n \geq 1}$ be a positive sequence which tends to $0$, and recall that 
\begin{align*}
&\mathcal{E}_{i}\cap{\hat{\mathcal{E}}_{i}}^c \\
& \quad =\left\{  V_{i}\leq
h_n\right\}  \cap\left\{  \hat{V}_{i}>h_n\right\}  \\
&  \quad \subset\left\{  V_{i}\leq h_n\right\}  \cap\left\{  h_n<V_{i}+\left\|
\left(  \hat{\Pi}_{D}-\Pi_{D}\right) \left(  \bx-\bX_{i}\right)  \right\|
\right\}  \\
&  \quad \subset\left\{V_{i}\leq h_n\right\}  \cap\left[  \left\{  h_n\left(
1-\kappa_n\right)  <V_{i}\right\}  \cup\left\{  \left\| \left(  \hat{\Pi
}_{D}-\Pi_{D}\right)  \left(  \bx-\bX_{i}\right)  \right\| >\kappa_nh_n\right\}
\right]  \\
&  \quad =\left\{  h_n\left(  1-\kappa_n\right)  <V_{i}\leq h_n\right\}  \\
& \quad \qquad \cup\left[
\left\{  V_{i}\leq h_n\right\}  \cap\left\{  \left\| \left(  \hat{\Pi
}_{D}-\Pi_{D}\right)  \left( \bx- \bX_{i}\right)  \right\| >\kappa_nh_n\right\}
\right].
\end{align*}
Thus%
\begin{align*}
& \sum_{i=1}^{n}{\mathbf 1}_{\mathcal{E}_{i}\cap{\hat
{\mathcal{E}}_{i}}^c}\\
& \quad \leq\sum_{i=1}^{n}{\mathbf 1}_{\left\{  h_n\left(
1-\kappa_n\right)  <V_{i}\leq h_n\right\}  }+\sum_{i=1}^{n}{\mathbf 1}_{\left\{
\left\| \left(  \hat{\Pi}_{D}-\Pi_{D}\right)  \left(  \bx-\bX_{i}\right)
\right\| >\kappa_nh_n\right\}  }{\mathbf 1}_{\left\{  V_{i}\leq h_n\right\}  },%
\end{align*}
and therefore
\begin{align}
\left[  \frac{1}{nF\left(  h_n\right)  }\sum_{i=1}^{n}\mathbf{1}%
_{\mathcal{E}_{i}\cap{\hat{\mathcal{E}}_{i}}^c}\right]  ^{2} &
\leq 2\left[  \frac{1}{nF\left(  h_n\right)  }\sum_{i=1}^{n}\mathbf{1}%
_{\left\{  h_n\left(  1-\kappa_n\right)  <V_{i}\leq h_n\right\}  }\right]  ^{2} \nonumber\\
&  \quad +2\left[  \sum_{i=1}^{n}{\mathbf 1}_{\left\{  \left\| \left(
\hat{\Pi}_{D}-\Pi_{D}\right)  \left( \bx-\bX_{i}\right)  \right\|
>\kappa_nh_n\right\}  }{\mathbf 1}_{\left\{  V_{i}\leq h_n\right\}  }\right]  ^{2}. \label{iwd}
\end{align}
Taking expectations and mimicking the method used in the proof
of Lemma \ref{L3}, we easily obtain
\begin{equation}
\mathbb{E}\left[  \frac{1}{nF\left(  h_n \right)  }\sum_{i=1}^{n}\mathbf{1}_{\left\{  h_n\left(  1-\kappa_n\right)  <V_{i}\leq h_n\right\}  }\right]  ^{2}\leq C\kappa_n^{2}.\label{B1}
\end{equation}
It remains to bound the last term on the right-hand side of (\ref{iwd}). To this aim, using the fact that $\bX$ has bounded support, we may write, for $\mu$-almost all $\bx$,
\begin{align*}
& {\mathbf 1}_{\left\{  \left\| \left(  \hat{\Pi}_{D}-\Pi_{D}\right)
\left( \bx- \bX_{i}\right)  \right\| >\kappa_nh_n\right\}  }{\mathbf 1}_{\left\{
V_{i}\leq h_n\right\}  }  \\
& \quad  \leq{\mathbf 1}_{\left\{  \left\| \hat{\Pi
}_{D}-\Pi_{D}\right\| _{\infty}>\kappa_n h_n/M\right\}  }{\mathbf 1}_{\left\{
V_{i}\leq h_n\right\}  }\\
& \qquad (\mbox{for some positive } M)\\
&\quad  \leq{\mathbf 1}_{\mathcal{A}_{n}^{c}}{\mathbf 1}_{\left\{  V_{i}\leq
h_n\right\}  }+{\mathbf 1}_{\left\{  \left\| \hat{\Pi}_{D}-\Pi
_{D}\right\| _{\infty}>\kappa_n h_n/M\right\}  \cap\mathcal{A}_{n}}%
{\mathbf 1}_{\left\{  V_{i}\leq h_n\right\}  },%
\end{align*}
where the set $\mathcal{A}_{n}$ is the same as in Lemma \ref{L0}. Clearly, by Cauchy-Schwarz inequality,
\begin{align}
& \mathbb{E}\left[  {\mathbf 1}_{\mathcal{A}_{n}^{c}}\left(  \frac
{1}{nF\left(  h_n\right)  }\sum_{i=1}^{n}{\mathbf 1}_{\left\{  V_{i}\leq
h_n\right\}  }\right) \right] ^2\nonumber  \\
& \quad  \leq\mathbb{P}^{1/2}\left(
\mathcal{A}_{n}^{c}\right)  \mathbb{E}^{1/2}\left[  \left(  \frac{1}{nF\left(
h_n\right)  }\sum_{i=1}^{n}{\mathbf 1}_{\left\{  V_{i}\leq h_n\right\}
}\right)  ^{4}\right] \nonumber \\
& \quad =\mathcal O \left(\exp\left(  -Cn\right) \right) \label{B2},
\end{align}
where the last inequality arises from Lemma \ref{L0}. It remains to bound the term
$$ \mathbb{E}\left[  \frac{1}{nF\left(  h_n\right)  }\sum_{i=1}%
^{n}{\mathbf 1}_{\left\{  \left\| \hat{\Pi}_{D}-\Pi_{D}\right\|
_{\infty}>\kappa_nh_n/M\right\}  \cap\mathcal{A}_{n}}{\mathbf 1}_{\left\{
V_{i}\leq h_n\right\}  } \right]^2.$$
We have
\begin{align*}
& \mathbb{E}\left[ \frac{1}{nF\left(  h_n\right)  }\sum_{i=1}%
^{n}{\mathbf 1}_{\left\{  \left\| \hat{\Pi}_{D}-\Pi_{D}\right\|
_{\infty}>\kappa_nh_n/M\right\}  \cap\mathcal{A}_{n}}{\mathbf 1}_{\left\{
V_{i}\leq h_n\right\}  }\right]^2  \\
& \quad=\mathbb{E}\left[  {\mathbf 1}_{\left\{  \left\| \hat{\Pi}%
_{D}-\Pi_{D}\right\| _{\infty}>\kappa_nh_n/M\right\}  \cap\mathcal{A}_{n}%
}\left(  \frac{1}{nF\left(  h_n\right)  }\sum_{i=1}^{n}{\mathbf 1}_{\left\{
V_{i}\leq h_n\right\}  }\right)  ^{2}\right].
\end{align*}
Using again Cauchy-Schwarz inequality and a bound on the fourth moment of $( 1/nF(  h_n))
\sum_{i=1}^{n}{\mathbf 1}_{\{  V_{i}\leq h_n\}  }$, it suffices to bound accurately %
\begin{align}
\mathbb{P}\left(  \left\{  \left\| \hat{\Pi}_{D}-\Pi_{D}\right\|
_{\infty}>\kappa_nh_n/M\right\}  \cap\mathcal{A}_{n}\right)    & \leq
\mathbb{P}\left(  \left\| \Gamma_{n}-\Gamma\right\| _{\infty}>\kappa_nh_n
C/M\right)  \nonumber\\
& = \mathcal O \left (\exp\left(  -C n\kappa_n^2h_n^{2}\right)\right),  \label{B3}%
\end{align}
where we used the bound (\ref{tango}) and statement $(i)$ in Theorem \ref{T1}.\\

Collecting the bounds (\ref{iwd})-(\ref{B1})-(\ref{B2})-(\ref{B3}), we finally obtain %
\[
\mathbb{E}\left[  \frac{1}{nF\left(  h_n\right)  }\sum_{i=1}^{n}\mathbf{1}_{\mathcal{E}_{i}\cap{\hat{\mathcal{E}}_{i}}^c}\right]
^{2}=\mathcal{O}\left(  \kappa_n^{2}\right)  +\mathcal{O}\left(  \exp\left(
-Cn\kappa_n^2h_n^{2}\right)  \right).
\]
The choice $\kappa_n^{\ast 2}\asymp\exp\left(  -Cn\kappa_n^{\ast 2}h_n^{2}\right)  $, i.e., 
$$\kappa_n^{\ast 2}\asymp \frac{\log\left(  nh^{2}\right)}{ nh_n^{2}}$$ 
leads to the desired result.\hfill $\blacksquare$

\end{document}